# HIGH ORDER MULTI-SCALE WALL-LAWS, PART I : THE PERIODIC CASE

DIDIER BRESCH† AND VUK MILISIC*

**Abstract.** In this work we present new wall-laws boundary conditions including microscopic oscillations. We consider a newtonian flow in domains with periodic rough boundaries that we simplify considering a Laplace operator with periodic inflow and outflow boundary conditions. Following the previous approaches, see [A. Mikelic, W. Jäger, J. Diff. Eqs, 170, 96–122, (2001) ] and [Y. Achdou, O. Pironneau, F. Valentin, J. Comput. Phys, 147, 1, 187–218, (1998)], we construct high order boundary layer approximations and rigorously justify their rates of convergence with respect to $\epsilon$ (the roughness' thickness). We establish mathematically a poor convergence rate for averaged second-order wall-laws as it was illustrated numerically for instance in [Y. Achdou, O. Pironneau, F. Valentin, J. Comput. Phys, 147, 1, 187–218, (1998)]. In comparison, we establish exponential error estimates in the case of explicit multi-scale ansatz. This motivates our study to derive implicit first order multi-scale wall-laws and to show that its rate of convergence is at least of order $\epsilon^{\frac{3}{2}}$. We provide a numerical assessment of the claims as well as a counter-example that evidences the impossibility of an averaged second order wall-law. Our paper may be seen as the first stone to derive efficient high order wall-laws boundary conditions.

**1. Introduction.** The main goal of wall-laws is to remove the stiff part from boundary layers, replacing the classical no-slip boundary condition by a more sophisticated relation between the variables and their derivatives. They are extensively used in numerical simulations to eliminate regions of strong gradients or regions of complex geometry (rough boundaries) from the domain of computation. Depending on the field of applications, (porous media, fluid mechanics, heat transfer, electromagnetism), wall-laws may be called BEAVERS-JOSEPH, SAFFMAN-JOSEPH, NAVIER, FOURIER, LEONTOVITCH type laws.

High order effective macroscopic boundary conditions may also be proposed if we choose a higher degree ansatz, see [7] for applications in microfluidic. In a similar perspective but in the context of fluid mechanics, numerical simulations have shown that second order macroscopic wall-laws provide the same order of approximation as the first order approximation. Recently a generalized wall-law formulation has been obtained for curved rough boundaries [17, 19] and for random roughness [4]. Note that such generalizations are important from a practical point of view when dealing with e.g. coastal effects in geophysical flows. From a mathematical point of view, wall-laws are also interesting. In the proof of convergence to the Euler equations, the 2D Navier-Stokes system is complemented with wall-laws of the NAVIER type [6]. Recently several papers analyze in various settings the properties of such boundary conditions, see [12], [16], [11], [5], [13].

In this paper, we focus on fluid flows. Starting from the Stokes system, we simplify the problem by studying the axial velocity through the resolution of a specific Poisson problem with periodic inlet and outlet boundary conditions. Our scope is to justify mathematically higher order macroscopic wall-laws and to explain why in their averaged form they do not provide better results than the first order laws. We shall explain how to get better estimates including some coefficients depending on the microscopic variables: this leads to new oscillating wall-laws.

†LAMA, UMR 5127 CNRS, Université de Savoie, 73217 Le Bourget du Lac cedex, FRANCE
*LJK-IMAG, UMR 5523 CNRS, 51 rue des Mathématiques, B.P.53, 38041 Grenoble cedex 9, FRANCE

The basic scheme to establish standard averaged wall-laws is the following (see fig. 1.1): First we use an ansatz for the velocity and the pressure which will give, after an adequate extension, a main order term completed with some boundary layer correctors defined on the whole rough domain (fig. 1.1, step I). This is possible due to the boundary layer theory that can be seen as a particular case of a general homogenization process. In a second time, a specific average is performed on this approximation and a new boundary condition of mixed type is recovered on a smooth fictitious interface strictly contained in the domain (fig. 1.1, step II). As one sees on the figure the only difference between Achdou's and Jäger's approaches is situated in the boundary layer's construction. It is an easy task to show that they are in fact a specific lift one of the other.

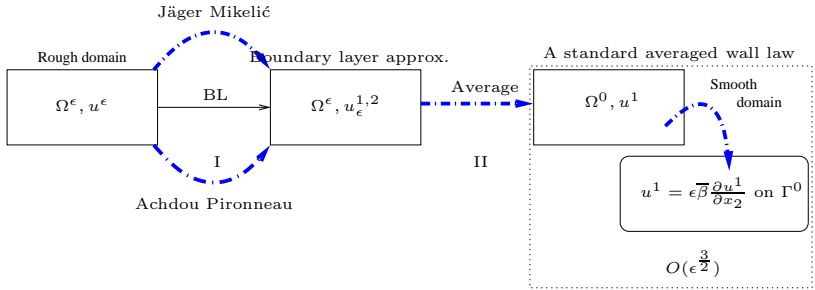

figure 1.1: The standard approach: from the exact solution to an averaged wall-law

The main result in our paper is the derivation of a high order boundary layer approximation that satisfies the homogeneous Dirichlet boundary condition on the rough wall and that leads to new wall-laws with microscopic effects see fig. 1.2. The ansatz is expanded up to the second order in $\epsilon$ and an exponential convergence in the interior domain is obtained using it, (fig. 1.2 step I'). Despite this great rate of convergence, the corresponding second order averaged wall-law behaves badly and does not preserve the nice convergence properties of full boundary layer approximations. The estimates show that this is due to the great influence of microscopic oscillations. We then derive new wall-laws that do converge exponentially on the smooth domain. They have the form of explicit non-homogeneous Dirichlet boundary conditions and they depend on the zeroth order Poiseuille flow as well as on the microscopic oscillations on the fictitious interface (fig. 1.2 step II').

At this stage, we go one step further and derive an implicit multi-scale first order wall-law. We obtain a SAFFMAN-JOSEPH's like law that now contains a coefficient that includes the microscopic oscillations. We rigorously derive a rate of convergence in $\epsilon^{\frac{3}{2}}$, thanks to the steps introduced in the previous sections (fig. 1.2 step II").

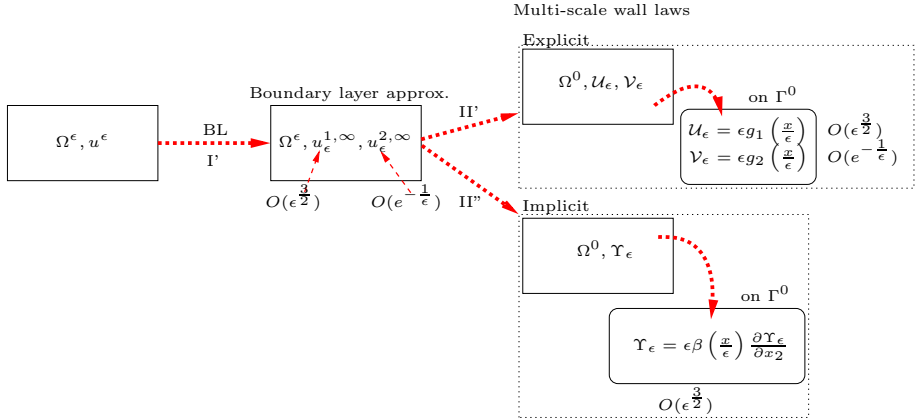

figure 1.2: The new approach: from the exact solution to multi-scale wall-laws

We underline that this work is a necessary building block when studying wall-laws for the stationary Navier-Stokes equations: asymptotic expansion of the quadratic non-linearity transfers a cascade of contributions to the microscopic cell problems, as already noticed in [2]. The first order cell problem is homogeneous and the second-order cell problem involves the non-linearity of the first order approximation. Until now, every averaged wall-law was only first order accurate and thus wall-laws were not able to display second order effects of non-linearities.

In a wider context that does not concern only fluid flows, the main concept this work emphasizes is the following: we have shown that it is possible to replace a geometrical roughness and "smooth" boundary conditions (in the sense unperturbed, as for instance homogeneous Dirichlet ones) by a smooth domain but with a multi-scale perturbed boundary conditions, (see fig. 1.3 below). Depending on the kind of boundary perturbation, we get different orders of precision in this process. For complex multi-scale 3D problems, we still expect some numerical gain when performing this switch, especially if one uses some increased multi-scale finite element bases (see [9] and references therein).

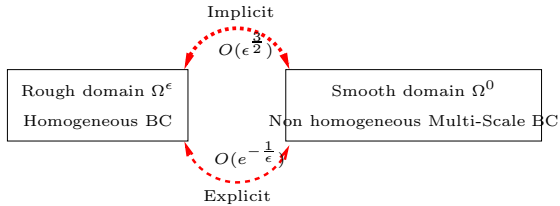

figure 1.3: One of the main points of this article: switching perturbations from geometry to boundary data. (BC stands for boundary conditions)

To show the practical importance of the results above, in Section 6, we perform numerical tests on a 2D case. For various values of $\epsilon$, we first compute the rough solution $u_\Delta^\epsilon$ on the whole domain $\Omega^\epsilon$, then we compute the wall-law solutions defined only on the interior smooth domain $\Omega^0$. We perform these tests in the periodic case. We recover exactly theoretical claims: numerical error estimates confirm that averaged wall-laws do not differ at first and second orders. We prove that our new implicit multi-scale wall-law provides better results than classical averaged laws. However, the fully explicit approximations still show higher order convergence rates with respect to $\epsilon$.

**2. The simplified problem: from Navier-Stokes to Laplace equation.** In this work, $\Omega^\epsilon$ denotes the rough domain in $\mathbb{R}^2$ depicted in fig. 2.1, $\Omega^0$ denotes the smooth one, $\Gamma^\epsilon$ is the rough boundary and $\Gamma^0$ (resp. $\Gamma^1$) the lower (resp. upper) smooth one (see fig. 2.1).

HYPOTHESES 2.1. *The rough boundary $\Gamma^\epsilon$ is described as a periodic repetition at the microscopic scale of a single boundary cell $P^0$. The latter can be parameterized as the graph of a Lipshitz function $f : [0, 2\pi[ \to [-1 : 0[$ such that*

$$P^0 = \{y \in [0, 2\pi] \times [-1 : 0[ \,/\, y_2 = f(y_1)\} \tag{2.1}$$

*Moreover we suppose that $f$ is negative definite, i.e. there exists a positive constant $\delta$ such that $f(y_1) < \delta$ for all $y_1 \in [0, 2\pi]$.* We assume that the ratio between $L$ (the width of $\Omega^0$) and $2\pi\epsilon$ (the width of the periodic cell) is always an integer called $N$. We

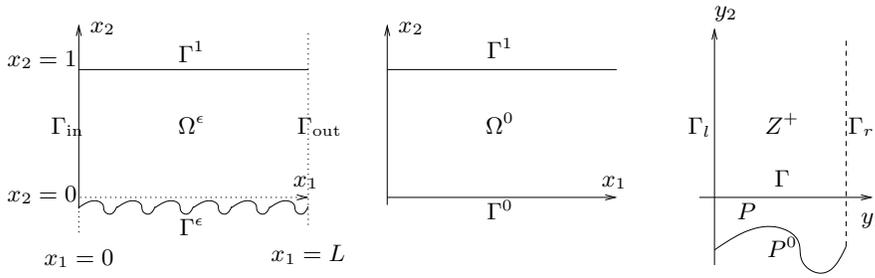

figure 2.1: *Rough, smooth and cell domains*

consider a simplified setting that avoids the theoretical difficulties and the non-linear complications of the full Navier-Stokes equations. Starting from the Stokes system, we consider a Poisson problem for the axial component of the velocity. The pressure gradient is assumed to reduce to a constant right hand side $C$. We consider only periodic inflow and outflow boundary conditions. The simplified formulation reads : find $u^\epsilon$ such that

$$\begin{cases} -\Delta u^\epsilon = C, & \text{for } x \in \Omega^\epsilon, \\ u^\epsilon = 0, & x \in \Gamma^\epsilon \cup \Gamma^1, \\ u^\epsilon \text{ is } x_1 \text{ periodic.} \end{cases} \tag{2.2}$$

We underline that the results below can be directly extended to rough domains with smooth holes and to the Stokes system.
In what follows, functions that do depend on $y = x/\epsilon$ should be indexed by an $\epsilon$ (e.g. $\mathcal{U}_\epsilon = \mathcal{U}_\epsilon(x, x/\epsilon)$).

**3. The full boundary layers correctors.**

**3.1. A zeroth order approximation.** When $\epsilon = 0$, the rough domain $\Omega^\epsilon$ reduces to $\Omega^0$ which is smooth. The solution of system (2.2) in this limit is known and explicit: it is the so-called Poiseuille profile :

$$\tilde{u}^0(x) = \frac{C}{2}(1 - x_2)x_2, \quad \forall x \in \Omega^0,$$

the latter term should be our zeroth order approximation when performing an asymptotic expansion w.r.t. $\epsilon$ for $\epsilon > 0$. The determining step is then how to extend this zeroth order approximation so that it is defined on the whole domain $\Omega^\epsilon$. A possible choice is to use the Taylor expansion of $\tilde{u}^0$ near $x_2 = 0$, it leads to define the zeroth order expansion as a $C^1(\Omega^\epsilon)$ function that reads

$$u_1^0(x) = \begin{cases} \tilde{u}^0(x), & \text{if } x \in \Omega^0 \\ \dfrac{\partial \tilde{u}^0}{\partial x_2}(x_1, 0)x_2, & \text{if } x \in \Omega^\epsilon \setminus \Omega^0. \end{cases}$$

Remark that this particular choice does not satisfy the homogeneous Dirichlet boundary condition on $\Gamma^\epsilon$. Next we estimate the zeroth order error wrt the exact solution.

PROPOSITION 1. *If $\Omega^\epsilon$ is a open connected piecewise smooth domain, the solution $u^\epsilon$ exists in $H^1(\Omega^\epsilon)$ and is unique. Moreover we have*

$$\left\| u^\epsilon - u_1^0 \right\|_{H^1(\Omega^\epsilon)} \leq c_1 \sqrt{\epsilon}, \quad \left\| u^\epsilon - u_1^0 \right\|_{L^2(\Omega^0)} \leq c_2 \epsilon,$$

*where the constants $c_1$ and $c_2$ are independent on $\epsilon$*

*Proof.* It is based on standard *a priori* estimates and a duality argument. The existence and uniqueness of $u^\epsilon$ are standard and left to the reader. We focus only on the error estimates. Namely, $r^0 := u^\epsilon - u^0$ satisfies

$$\begin{cases} -\Delta r^0 = C\chi_{[\Omega^\epsilon \setminus \Omega^0]}, & \text{in } \Omega^\epsilon, \\ r^0 = 0, & \text{on } \Gamma^1, \\ r^0 = -\dfrac{\partial \tilde{u}^0}{\partial x_2}(x_1, 0)x_2 & \text{on } \Gamma^0, \\ r^0 \text{ is } x_1 - \text{periodic on } \Gamma_{\text{in}} \cup \Gamma_{\text{out}} \end{cases}$$

There, one remarks that a part of the error comes from the source term localized in $\Omega^\epsilon \setminus \Omega^0$, and another part comes from the non homogeneous boundary term on $\Gamma^\epsilon$. We set the lift

$$s = -\frac{\partial \tilde{u}^0}{\partial x_2} x_2 \chi_{[\Omega^\epsilon \setminus \Omega^0]}, \text{ and } z := r^0 - s,$$

then the weak formulation reads :

$$(\nabla z, \nabla v)_{\Omega^\epsilon} = (C, v)_{\Omega^\epsilon \setminus \Omega^0} - (C, v)_{\Gamma^0}, \quad v \in H_0^1(\Omega^\epsilon),$$

where the last term in the rhs comes when applying the Laplace operator $\Delta$ on $s$. Thanks to Poincaré-like estimates we have the following properties of the $L^2$ norm and the $H^1$ semi-norm on $\Omega^\epsilon \setminus \Omega^0$

$$\left| (C, v)_{\Omega^\epsilon \setminus \Omega^0} - (C, v)_{\Gamma^0} \right| \leq c_3 \sqrt{\epsilon} \left( \int_{\Omega^\epsilon \setminus \Omega^0} v^2 \right)^{\frac{1}{2}} + c_4 \left( \int_{\Gamma^0} v^2 \right)^{\frac{1}{2}}$$

$$\leq c_5(\epsilon + \sqrt{\epsilon}) \left( \int_{\Omega^\epsilon \setminus \Omega^0} |\nabla v|^2 \right)^{\frac{1}{2}} \leq 2c_5 \sqrt{\epsilon} \|v\|_{H^1(\Omega^\epsilon)}.$$

This leads to the $H^1(\Omega^\epsilon)$ estimate. For the $L^2$ norm, we use the concept of a very weak solution [18]. Namely, one solves the dual problem: for a given $\varphi \in L^2(\Omega^0)$, $\varphi$ being $x_1$ − periodic on $\Gamma_{\text{in}} \cup \Gamma_{\text{out}}$ find $v \in H^2(\Omega^0)$ such that

$$\begin{cases} -\Delta v = \varphi, & \forall x \in \Omega^0, \\ v = 0, & \forall x \in \Gamma^0 \cup \Gamma^1, \\ v \text{ is } x_1 - \text{periodic on } \Gamma_{\text{in}} \cup \Gamma_{\text{out}}. \end{cases}$$

Considering the $L^2(\Omega^0)$ scalar product, and using the Green formula

$$\begin{aligned}(\varphi, r^0)_{\Omega^0} &= -(\Delta v, r^0)_{\Omega^0} = \left\langle \frac{\partial r^0}{\partial \mathbf{n}}, v \right\rangle_{\partial\Omega^0} - \left( \frac{\partial v}{\partial \mathbf{n}}, r^0 \right)_{\partial\Omega^0} - (v, \Delta r^0)_{\Omega^0}, \\ &= \left\langle v, \frac{\partial r^0}{\partial \mathbf{n}} \right\rangle_{\Gamma_{\text{in}} \cup \Gamma_{\text{out}}} - \left( \frac{\partial v}{\partial \mathbf{n}}, r^0 \right)_{\Gamma^0 \cup \Gamma^1}, \end{aligned} \quad (3.1)$$

where the brackets refer to the dual product in $H^{\frac{1}{2}}(\Gamma^0)$, and the rest of products are in $L^2$, either on $\Gamma^0$ or on $\Omega^0$. Then, one computes

$$|(\varphi, r^0)| \leq \left| \left( \frac{\partial v}{\partial \mathbf{n}}, r^0 \right)_{\Gamma^0} \right| \leq c_6 \|\varphi\|_{L^2(\Omega^0)} \|r^0\|_{L^2(\Gamma^0)}.$$

The last estimate is obtained thanks to a linear dependence of the normal derivative of the trace of $v$ on the data $\varphi$, [18]. Thanks to Poincaré estimates, one writes

$$\|r^0\|_{L^2(\Gamma^0)} \leq c_7 \sqrt{\epsilon} \|r^0\|_{H^1(\Omega^\epsilon \setminus \Omega^0)} \leq c_8 \sqrt{\epsilon} \|r^0\|_{H^1(\Omega^\epsilon)}$$

which ends the proof by taking the sup over all $\varphi$ in $L^2(\Omega^0)$. □

**3.2. A first order correction.** The zeroth order correction contains two distinct sources of errors : a part is due to the order of the extension in $\Omega^\epsilon \setminus \Omega^0$ and another part comes from a non homogeneous rest on $\Gamma^\epsilon$. In what follows we show that a first order extension $u_1^0$ can be corrected by series of terms that makes the full boundary layer approximation vanish on $\Gamma^\epsilon$.

*The micrscopic cell problem* : In order to correct $u_1^0$ on $\Gamma^\epsilon$, one starts by solving a microscopic cell problem that reads : find $\beta$ s.t.

$$\begin{cases} -\Delta \beta = 0, & \text{in } Z^+ \cup P, \\ \beta = -y_2, & \text{on } P^0, \\ \beta \text{ is } y_1 - \text{periodic}. \end{cases} \quad (3.2)$$

We define the microscopic average along the fictitious interface $\Gamma$ :

$$\overline{\beta} = \frac{1}{2\pi} \int_0^{2\pi} \beta(y_1, 0) dy_1.$$

As $Z^+ \cup P$ is unbounded in the $y_2$ direction, we define

$$D^{1,2} = \{v \in L^1_{\text{loc}}(Z^+ \cup P) / Dv \in L^2(Z^+ \cup P)^2, v \text{ is } y_1 - \text{periodic }\},$$

then one has the following result :

THEOREM 3.1. *Under hypotheses 2.1, there exists $\beta$, a unique solution of (3.2) belonging to $D^{1,2}$. Moreover, there exists a unique periodic solution $\eta \in H^{\frac{1}{2}}(\Gamma)$, of the following problem*

$$< S\eta, \mu > = < 1, \mu >, \quad \forall \mu \in H^{\frac{1}{2}}(\Gamma),$$

where $<,>$ is the $H^{-\frac{1}{2}}(\Gamma) - H^{\frac{1}{2}}(\Gamma)$ duality bracket, and $S$ the inverse of the Steklov-Poincaré operator (see appendix A.1). One has the following correspondance between $\beta$ and the interface solution $\eta$ :

$$\beta = H_{Z^+}\eta + H_P\eta,$$

where $H_{Z^+}\eta$ (resp. $H_P\eta$) is the $y_1$-periodic harmonic extension of $\eta$ on $Z^+$ (resp. $P$). The solution in $Z^+$ can be written explicitly as a series of Fourier coefficients of $\eta$ and reads :

$$H_{Z^+}\eta = \beta(y) = \sum_{k=-\infty}^{\infty} \eta_k e^{iky_1 - |k|y_2}, \quad \forall y \in Z^+, \quad \eta_k = \int_0^{2\pi} \eta(y_1) e^{-iky_1} dy_1.$$

In the macroscopic domain $\Omega^0$ this leads to

$$\left\| \beta\left(\frac{\cdot}{\epsilon}\right) - \overline{\beta} \right\|_{L^2(\Omega^0)} \leq K\sqrt{\epsilon} \|\eta\|_{H^{\frac{1}{2}}(\Gamma)}.$$

The proof is given in the appendix for sake of conciseness. The corresponding macroscopic full boundary layer corrector should contain at this stage

$$u_1^0 + \epsilon \frac{\partial u_1^0}{\partial x_2}(x_1, 0) \left(\beta\left(\frac{x}{\epsilon}\right) - \overline{\beta}\right),$$

where we subtract $\overline{\beta}$ in order to cancel $\beta$'s errors on $\Gamma^1$.
In order to cancel the contribution of the constant $\overline{\beta}$ near the rough boundary but keep its benefit close to $\Gamma^1$, one solves the "counter-flow" problem: find $d$ s.t.

$$\begin{cases} -\Delta d = 0, & \text{in } \Omega^0, \\ d = 1 \text{ on } \Gamma^0, d = 0 \text{ on } \Gamma^1, \\ d \text{ is } x_1 - \text{periodic on } \Gamma_{\text{in}} \cup \Gamma_{\text{out}}, \end{cases} \quad (3.3)$$

the solution is explicit and reads $d = (1 - x_2)$. Moreover, it can be extended to the whole domain $\Omega^\epsilon$. The complete first order approximation now reads :

$$u_\epsilon^{1,2} := u_1^0 + \epsilon \frac{\partial u_1^0}{\partial x_2}(x_1, 0)(\beta - \overline{\beta}) + \epsilon \frac{\partial u_1^0}{\partial x_2}(x_1, 0)\overline{\beta}(1 - x_2), \quad \forall x \in \Omega^\epsilon$$

$$= u_1^0 + \epsilon \frac{\partial u_1^0}{\partial x_2}(x_1, 0)(\beta - \overline{\beta} x_2),$$

the first index of $u_\epsilon^{1,2}$ corresponds to the extension order of $\tilde{u}^0$ in $\Omega^\epsilon \setminus \Omega^0$, while the second index is the order of the error on $\Gamma^\epsilon$. Indeed, if we consider the trace of $u_\epsilon^{1,2}$ on $\Gamma^\epsilon$, we have a second order error

$$u_\epsilon^{1,2}\big|_{\Gamma^\epsilon} = \epsilon^2 \left(\frac{\partial u_1^0}{\partial x_2}\overline{\beta}\right) \frac{x_2}{\epsilon} = \epsilon^2 \left(\frac{\partial u_1^0}{\partial x_2}\overline{\beta}\right) y_2.$$

Again, this error is linear and should be corrected by the micro boundary layer $\beta$. A similar macroscopic boundary layer correction process should be performed at any order leading to

$$u_\epsilon^{1,\infty} = u_1^0 + \epsilon \frac{\partial u_1^0}{\partial x_2}(x_1, 0) \left[\left(\beta\left(\frac{x}{\epsilon}\right) - \overline{\beta} x_2\right) + \epsilon\overline{\beta}\left(\beta\left(\frac{x}{\epsilon}\right) - \overline{\beta} x_2\right)\right.$$

$$\left. - \epsilon^2 \overline{\beta}^2 \left(\beta\left(\frac{x}{\epsilon}\right) - \overline{\beta} x_2\right) + \ldots\right] \quad (3.4)$$

$$= u_1^0 + \frac{\epsilon}{1 + \epsilon\overline{\beta}} \frac{\partial u_1^0}{\partial x_2}(x_1, 0) \left(\beta\left(\frac{x}{\epsilon}\right) - \overline{\beta} x_2\right).$$

This approximation satisfies a homogeneous Dirichlet boundary condition on $\Gamma^\epsilon$, and solves

$$\begin{cases} -\Delta u_\epsilon^{1,\infty} = C\chi_{[\Omega^\epsilon \setminus \Omega^0]}, & \text{in } \Omega^\epsilon, \\ u_\epsilon^{1,\infty} = 0, & \text{on } \Gamma^\epsilon, \\ u_\epsilon^{1,\infty} = \dfrac{\epsilon}{1+\epsilon\overline{\beta}}\dfrac{\partial u_1^0}{\partial x_2}\left(\beta\left(\dfrac{x_1}{\epsilon},0\right) - \overline{\beta}\right), & \text{on } \Gamma^\epsilon, \\ u_\epsilon^{1,\infty} \text{ is } x_1-\text{periodic on } \Gamma_{\text{in}} \cup \Gamma_{\text{out}}. \end{cases} \quad (3.5)$$

If we consider the corresponding approximation error, we obtain

PROPOSITION 2. *Under hypotheses 2.1, the error of the first order approximation satisfies*

$$\left\|u^\epsilon - u_\epsilon^{1,\infty}\right\|_{H^1(\Omega^\epsilon)} \le c_8\epsilon, \quad \left\|u^\epsilon - u_\epsilon^{1,\infty}\right\|_{L^2(\Omega^0)} \le c_9\epsilon^{\frac{3}{2}},$$

*where the constants $c_8, c_9$ are independent on $\epsilon$.*

The proof follows the same lines as in proposition 1 except that the significant source of errors is the rhs of the first equation in (3.5), while an exponentially small microscopic perturbation lies on $\Gamma^1$, on the contrary there are no errors on $\Gamma^\epsilon$, because $u_\epsilon^{1,\infty}=0$ on it.

**3.3. Second order approximation.** Instead of extending only linearly the Poiseuille profile it is obvious that a quadratic term is missing to complete the approximation. In the following $u_2^0$ denotes the second order extension of $\tilde{u}^0$ in $\Omega^\epsilon \setminus \Omega^0$.

$$u_2^0 := \begin{cases} \tilde{u}^0, & x \in \Omega^0 \\ \dfrac{\partial \tilde{u}^0}{\partial x_2}(x_1,0)x_2 + \dfrac{\partial^2 \tilde{u}^0}{\partial x_2^2}(x_1,0)\dfrac{x_2^2}{2}, & x \in \chi_{[\Omega^\epsilon\setminus\Omega^0]} \end{cases} = \dfrac{C}{2}(1-x_2)x_2, \quad \forall x \in \Omega^\epsilon.$$

The second order error on $\Gamma^\epsilon$ is corrected thanks to a new cell problem : find $\gamma \in D^{1,2}$ solving

$$\begin{cases} -\Delta \gamma = 0, & \text{in } Z^+ \cup P, \\ \gamma = -y_2^2, & \text{on } P^0, \\ \gamma \text{ periodic in } y_1. \end{cases} \quad (3.6)$$

The proof of the following proposition is left in the appendix A.2.

PROPOSITION 3. *Under hypotheses 2.1, there exists a unique solution $\gamma$ of (3.6) in $D^{1,2}(Z^+ \cup P)$. Moreover it admits a power series of Fourier modes in $Z^+$ and $\gamma \in [-1,0]$ if $P \subset [0,2\pi] \times [-1,0]$.*

The horizontal average is denoted $\overline{\gamma}$. The same multi-scale process leads to write the full boundary layer approximation as

$$u_\epsilon^{2,3} = u_2^0 + \dfrac{\epsilon}{1+\epsilon\overline{\beta}}\dfrac{\partial u_2^0}{\partial x_2}(x_1,0)\left(\beta\left(\dfrac{x}{\epsilon}\right) - \overline{\beta}x_2\right) + \dfrac{\epsilon^2}{2}\dfrac{\partial^2 u_2^0}{\partial x_2}(x_1,0)\left(\gamma\left(\dfrac{x}{\epsilon}\right) - \overline{\gamma}x_2\right).$$

Again a third error remains on $\Gamma^\epsilon$ and it is linear wrt to $y_2$, thus it should be corrected thanks to the series of first order cell problems as in (3.4). We set $u_\epsilon^{2,\infty}$ to be the second order approximation that satisfies a homogeneous Dirichlet boundary condition on $\Gamma^\epsilon$,

it reads :

$$u_\epsilon^{2,\infty} = u_2^0 + \frac{\epsilon}{1+\epsilon\overline{\beta}}\frac{\partial u_2^0}{\partial x_2}(x_1,0)\left(\beta\left(\frac{x}{\epsilon}\right) - \overline{\beta}x_2\right)$$
$$+ \frac{\epsilon^2}{2}\frac{\partial^2 u_2^0}{\partial x_2}(x_1,0)\left[\left(\gamma\left(\frac{x}{\epsilon}\right) - \overline{\gamma}x_2\right) + \frac{\epsilon\overline{\gamma}}{1+\epsilon\overline{\beta}}\left(\beta\left(\frac{x}{\epsilon}\right) - \overline{\beta}x_2\right)\right].$$

Our approximation satisfies the following boundary value problem

$$\begin{cases} -\Delta u_\epsilon^{2,\infty} = C, \text{ in } \Omega^\epsilon, \\ u_\epsilon^{2,\infty} = 0, \text{ on } \Gamma^\epsilon, \\ u_\epsilon^{2,\infty} = g_\epsilon, \text{ on } \Gamma^1, \\ u_\epsilon^{2,\infty} \text{ is } x_1 - \text{periodic on } \Gamma_{\text{in}} \cup \Gamma_{\text{out}}, \end{cases} \quad (3.7)$$

where $g$ is the contribution of the microscopic correctors on $\Gamma^1$ and reads :

$$g_\epsilon = \frac{\partial u_2^0}{\partial x_2}(x_1,0)\left(\beta\left(\frac{x_1}{\epsilon},1\right) - \overline{\beta}\right)$$
$$+ \frac{\epsilon^2}{2}\frac{\partial^2 u_2^0}{\partial x_2}(x_1,0)\left[\left(\gamma\left(\frac{x_1}{\epsilon},1\right) - \overline{\gamma}\right) + \frac{\epsilon\overline{\gamma}}{1+\epsilon\overline{\beta}}\left(\beta\left(\frac{x_1}{\epsilon},1\right) - \overline{\beta}\right)\right].$$

Remark that the only error remains on $\Gamma^1$ and as the proposition below claims, it is exponentially small wrt $\epsilon$.

PROPOSITION 4. *Under hypotheses 2.1 the error of the first second order approximation satisfies*

$$\left\|u^\epsilon - u_\epsilon^{2,\infty}\right\|_{H^1(\Omega^\epsilon)} \leq c_{10}e^{-\frac{1}{\epsilon}}, \quad \left\|u^\epsilon - u_\epsilon^{2,\infty}\right\|_{L^2(\Omega^0)} \leq c_{11}\sqrt{\epsilon}e^{-\frac{1}{\epsilon}}.$$

*where the constants $c_6, c_7$ are independent on $\epsilon$.* The proof is identical to the one of proposition 1 except that the only source of errors is the contribution of function $g_\epsilon$, there are nor errors on $\Gamma^\epsilon$, neither source terms inside $\Omega^\epsilon$.

**4. Averaged wall-laws.**

**4.1. The averaged wall-laws: a new derivation process.** At this stage, we rewrite our first and second order approximations separating slow and fast variables

$$u_\epsilon^{1,\infty} = u_1^0 + \frac{\epsilon\overline{\beta}}{1+\epsilon\overline{\beta}}\frac{\partial u_1^0}{\partial x_2}(x_1,0)(1-x_2) + \frac{\epsilon}{1+\epsilon\overline{\beta}}\frac{\partial u_1^0}{\partial x_2}(x_1,0)\left(\beta\left(\frac{x}{\epsilon}\right) - \overline{\beta}\right),$$

$$u_\epsilon^{2,\infty} = u_2^0 + \frac{\epsilon\overline{\beta}}{1+\epsilon\overline{\beta}}\frac{\partial u_2^0}{\partial x_2}(x_1,0)(1-x_2)$$
$$+ \frac{\epsilon^2}{2}\frac{\partial^2 u_2^0}{\partial x_2}(x_1,0)\left[\overline{\gamma}(1-x_2) + \frac{\epsilon\overline{\gamma}\overline{\beta}}{1+\epsilon\overline{\beta}}(1-x_2)\right]$$
$$+ \frac{\epsilon}{1+\epsilon\overline{\beta}}\frac{\partial u_2^0}{\partial x_2}(x_1,0)\left(\beta\left(\frac{x}{\epsilon}\right) - \overline{\beta}\right)$$
$$+ \frac{\epsilon^2}{2}\frac{\partial^2 u_2^0}{\partial x_2}(x_1,0)\left[\left(\gamma\left(\frac{x}{\epsilon}\right) - \overline{\gamma}x_2\right) + \frac{\epsilon\overline{\gamma}}{1+\epsilon\overline{\beta}}\left(\beta\left(\frac{x}{\epsilon}\right) - \overline{\beta}x_2\right)\right].$$

We define the average wrt the fast variable in the horizontal direction:

$$\overline{v}(x) = \frac{1}{2\pi\epsilon}\int_0^{2\pi\epsilon} v(x_1+y, x_2)dy, \quad \forall v \in H^1(\Omega^\epsilon).$$

Then, one can see easily that for any $x$ in $\Omega^0$

$$\overline{u_\epsilon^{1,\infty}} = u_1^0 + \frac{\epsilon\overline{\beta}}{1+\epsilon\overline{\beta}}\frac{\partial u_1^0}{\partial x_2}(x_1,0)(1-x_2) =: u^1,$$

$$\overline{u_\epsilon^{2,\infty}} = u_2^0 + \frac{\epsilon\overline{\beta}}{1+\epsilon\overline{\beta}}\frac{\partial u_2^0}{\partial x_2}(x_1,0)(1-x_2)$$
$$+ \frac{\epsilon^2}{2}\frac{\partial^2 u_2^0}{\partial x_2}(x_1,0)\left[\overline{\gamma}(1-x_2) + \frac{\epsilon\overline{\gamma}\overline{\beta}}{1+\epsilon\overline{\beta}}(1-x_2)\right] =: u^2.$$

This means that the averaging process cancels the oscilations providing only macroscopic terms still depending on $\epsilon$. Moreover one has the following compact form of the full boundary layer correctors

$$
\begin{aligned}
u_\epsilon^{1,\infty} &= u^1 + \epsilon\frac{\partial u^1}{\partial x_2}(x_1,0)\left(\beta\left(\frac{x}{\epsilon}\right) - \overline{\beta}\right) \\
u_\epsilon^{2,\infty} &= u^2 + \epsilon\frac{\partial u^2}{\partial x_2}(x_1,0)\left(\beta\left(\frac{x}{\epsilon}\right) - \overline{\beta}\right) + \frac{\epsilon^2}{2}\frac{\partial^2 u^2}{\partial x_2^2}(x_1,0)\left(\gamma\left(\frac{x}{\epsilon}\right) - \overline{\gamma}\right).
\end{aligned}
\quad (4.1)
$$

At this point, if one computes the boundary value problem that $u^1$ and $u^2$ solve in the smooth domain, we obtain the two following Robin and Wentzel type problems. Namely, $u^1$ solves:

$$
\begin{cases}
-\Delta u^1 = C, & \forall x \in \Omega^0, \\
u^1 = \epsilon\overline{\beta}\dfrac{\partial u^1}{\partial x_2}, & \forall x \in \Gamma^0, \quad u^1 = 0, \quad \forall x \in \Gamma^1, \\
u^1 \text{ is } x_1-\text{periodic on } \Gamma_{\text{in}} \cup \Gamma_{\text{out}},
\end{cases}
\quad (4.2)
$$

whose explicit solution reads:

$$u^1(x) = -\frac{C}{2}\left(x_2^2 - \frac{x_2}{1+\epsilon\overline{\beta}} - \frac{\epsilon\overline{\beta}}{1+\epsilon\overline{\beta}}\right), \quad (4.3)$$

while the second order wall-law $u^2$ satisfies the folowing boundary value problem

$$
\begin{cases}
-\Delta u^2 = C, & \forall x \in \Omega^0, \\
u^2 = \epsilon\overline{\beta}\dfrac{\partial u^2}{\partial x_2} + \dfrac{\epsilon^2}{2}\overline{\gamma}\dfrac{\partial^2 u^2}{\partial x_2^2}, & \forall x \in \Gamma^0, \\
u^2 = 0, \quad \forall x \in \Gamma^1, u^2 \text{ is } x_1-\text{periodic on } \Gamma_{\text{in}} \cup \Gamma_{\text{out}}.
\end{cases}
\quad (4.4)
$$

**4.2. Existence and uniqueness of the second order wall-law.** Because problem (4.4) contains second order normal derivatives as components of the boundary condition, (in the literature this kind of boundary conditions are called of Wentzell boundary conditions) the existence and uniqueness is not a standard result. Here we provide it. First we transform the second-order normal boundary term in a tangential term of the same order. Then using the appropriate test function space, we can apply Green's formula on tangential directions and symmetrise the bilinear form associated to the problem.

LEMMA 4.1. *Under hypotheses 2.1,the system (4.4) admits a unique solution in* $H_\#^{1,1}(\Omega^0) = \{v \in H_{\Gamma^1}^1(\Omega^0); v \in H^1(\Gamma^0)\}$, *where* $H_{\Gamma^1}^1$ *is the set of functions belonging to* $H^1(\Omega^0)$, $x_1-$ *periodic on* $\Gamma_{\text{in}} \cup \Gamma_{\text{out}}$ *and vanishing on* $\Gamma^1$.

*Proof.* The boundary condition shall be transformed thanks to the first equation of (4.4) into

$$u = \epsilon\overline{\beta}\frac{\partial u}{\partial x_2} + \frac{\epsilon^2}{2}\overline{\gamma}\frac{\partial^2 u}{\partial x_2^2} = \epsilon\overline{\beta}\frac{\partial u}{\partial x_2} + \frac{\epsilon^2}{2}\overline{\gamma}\left(-C - \frac{\partial^2 u}{\partial x_1^2}\right), \quad \forall x \in \Gamma^0.$$

Because $P^0$ does not intersect $\Gamma$, and thanks to the maximum principle, $\beta > 0$ a.e. in $Z^+ \cup P$. This implies that $\overline{\beta} > 0$ which allows the weak formulatiuon [10] :

$$\frac{1}{\epsilon\overline{\beta}}(u,v)_{\Gamma^0} + (\nabla u, \nabla v)_{\Omega^0} - \epsilon\frac{\overline{\gamma}}{2\overline{\beta}}\left[\left(\frac{\partial u}{\partial x_1}v\right)(x_1,0)\right]_{x_1=0}^{x_1=L} - \epsilon\frac{\overline{\gamma}}{2\overline{\beta}}\left(\frac{\partial u}{\partial x_1}\frac{\partial v}{\partial x_1}\right)_{\Gamma^0}$$
$$= (C,v)_{\Omega^0} - \epsilon\frac{\overline{\gamma}}{2\overline{\beta}}(C,v)_{\Gamma^0},$$

where the third term of the lhs vanishes thanks to the periodicity of the solution and of the corresponding test functions of $H^1_\#(\Omega^0)$. We have obtained a symmetric problem. Because $\overline{\gamma} \in [-1,0[$ and $\overline{\beta} \in ]0,1]$, setting

$$a(u,v) = \frac{1}{\epsilon\overline{\beta}}(u,v)_{\Gamma^0} + (\nabla u, \nabla v)_{\Omega^0} - \epsilon\frac{\overline{\gamma}}{2\overline{\beta}}\left(\frac{\partial u}{\partial x_1}\frac{\partial v}{\partial x_1}\right)_{\Gamma^0}, v \in H^{1,1}_\#(\Omega^0),$$
$$l(v) = (C,v)_\Omega - \epsilon\frac{\overline{\gamma}}{2\overline{\beta}}(C,v)_{\Gamma^0},$$

one obtains a variational formulation where $a$ is coercive, $H^{1,1}_\#(\Omega^0)$ being endowed with the norm :

$$\|u\|_{H^{1,1}_\#(\Omega^0)} = \|u\|_{H^1(\Omega^0)} + \|u\|_{H^1(\Gamma^0)}.$$

Moreover, $a$ and $l$ are continuous on $H^{1,1}_\#(\Omega^0)$, thus the problem is solvable by the Lax-Milgram theorem. By the way, we derive the following energy estimates that describe the dependence of various norms upon $\epsilon$ :

$$\|u\|_{L^2(\Gamma^0)} \leq \sqrt{\epsilon}C, \quad \left\|\frac{\partial u}{\partial x_1}\right\|_{L^2(\Gamma^0)} \leq \frac{C}{\sqrt{\epsilon}}.$$

Note that when $\epsilon$ goes to zero, our approximation leaves $H^{1,1}_\#(\Omega^0)$ moving to $H^1_{\Gamma^1 \cup \Gamma^0}(\Omega^0)$: we loose the control over the tangential derivative on the boundary. □

In the particular case of a straight domain $\Omega^0$ this unique solution is explicit and reads

$$u^2(x) = -\frac{C}{2}\left(x_2^2 - \frac{x_2(1+\epsilon^2\overline{\gamma})}{1+\epsilon\overline{\beta}} - \frac{\epsilon(\overline{\beta}-\epsilon\overline{\gamma})}{1+\epsilon\overline{\beta}}\right). \tag{4.5}$$

**4.3. Macroscopic error estimate.** When replacing the Poiseuille profile in $\Omega^0$ by $u^1$ or $u^2$, one can compute the corresponding error estimates.

PROPOSITION 5. *Let $u^\epsilon$ be the solution of (2.2) and $u^1$ (resp. $u^2$) be the solution of (4.2) (resp. (4.4)), then Under hypotheses 2.1,*

$$\left\|u^\epsilon - u^1\right\|_{L^2(\Omega^0)} \leq C\epsilon^{\frac{3}{2}}, \text{ and } \left\|u^\epsilon - u^2\right\|_{L^2(\Omega^0)} \leq C\epsilon^{\frac{3}{2}}.$$

*Proof.* We only compute the error of the second order approximation, the case of $u^1$ is identical. We take advantage of estimates obtained in proposition 2 by inserting the full boundary layer corrector $u_\epsilon^{2,\infty}$ between $u^\epsilon$ and $u^2$ :

$$u^\epsilon - u^2 = u^\epsilon - u_\epsilon^{2,\infty} + u_\epsilon^{2,\infty} - u^2$$

$$= u^\epsilon - u_\epsilon^{2,\infty} + \epsilon \frac{\partial u^2}{\partial x_2}(x_1, 0)\left(\beta\left(\frac{x}{\epsilon}\right) - \overline{\beta}\right) + \frac{\epsilon^2}{2}\frac{\partial^2 u^2}{\partial x_2^2}(x_1, 0)\left(\gamma\left(\frac{x}{\epsilon}\right) - \overline{\gamma}\right),$$

where we used the compact form exhibited in (4.1). Then, one gets

$$\|u^\epsilon - u^2\|_{L^2(\Omega^0)} \leq \|u^\epsilon - u_\epsilon^{2,\infty}\|_{L^2(\Omega^0)} + K\epsilon\left((1+\epsilon^2)\|\beta - \overline{\beta}\|_{L^2(\Omega^0)} + \epsilon\|\gamma - \overline{\gamma}\|_{L^2(\Omega^0)}\right).$$

Thanks to proposition 2, and the last estimate in the claim of theorem 3.1, one gets the desired result. □

REMARK 4.1. *This result is crucial: it shows that the oscillations of the first order boundary layer $\epsilon \partial u^0 / \partial x_2 (\beta - \overline{\beta})$ are larger than the second order macroscopic contribution. It is also optimal (see section 6 for a numerical evidence). This observation motivates the sections below.*

**5. Multi-scale wall-laws.** In this section we continue the investigation in the sense introduced above. We aim to compute a solution that exists in $\Omega^0$ as $u^1$ or $u^2$ but that performs a better approximation of the exact solution $u^\epsilon$ restricted to $\Omega^0$. Below we shall show that this concept provides some new multi-scale wall-laws.

**5.1. The first order explicit wall-law .** How can first order correction be improved if the non-oscillating second order extension of SAFFMAN-JOSEPH's condition does not help. The aswer below will be to take into account some multi-scale features. If we consider the full boundary layer corrector $u_\epsilon^{1,\infty}$, it solves (3.5). Moreover, on the fictitious boundary $\Gamma^0$, its value is easily computed, namely

$$u_\epsilon^{1,\infty}\big|_{x_2=0} = \left\{u^1 + \epsilon \frac{\partial u^1}{\partial x_2}(x_1, 0)\left(\beta - \overline{\beta}\right)\right\}\bigg|_{x_2=0} = \epsilon \frac{\partial u^1}{\partial x_2}(x_1, 0)\beta(x_1, 0).$$

We use this value as a non-homogenous Dirichlet boundary condition on $\Gamma^0$ for a Poisson problem that is nevertheless homogeneous on $\Gamma^1$. Indeed, we consider the following problem

$$\begin{cases} -\Delta \mathcal{U}_\epsilon = C, & \forall x \in \Omega^0, \\ \mathcal{U}_\epsilon = \epsilon \frac{\partial u^1}{\partial x_2}(x_1, 0)\beta\left(\frac{x_1}{\epsilon}, 0\right), & \forall x \in \Gamma^0, \\ \mathcal{U}_\epsilon = 0, & \forall x \in \Gamma^1, \quad \mathcal{U}_\epsilon \text{ is } x_1 - \text{periodic on } \Gamma_{\text{in}} \cup \Gamma_{\text{out}}, \end{cases} \quad (5.1)$$

and we claim the following

PROPOSITION 6. *Under hypotheses 2.1, one gets the following error estimates*

$$\|u^\epsilon - \mathcal{U}_\epsilon\|_{L^2(\Omega^0)} \leq c_{12}\epsilon^{\frac{3}{2}}.$$

*Proof.* Following the same lines as in the proof of proposition 5, one inserts the full boundary layer approximation error $r^{1,\infty} := u^\epsilon - u_\epsilon^{1,\infty}$ :

$$r_{\text{bl}}^1 = u^\epsilon - u_\epsilon^{1,\infty} + u_\epsilon^{1,\infty} - \mathcal{U}_\epsilon = r^{1,\infty} - \left[\mathcal{U}_\epsilon - u_\epsilon^{1,\infty}\right] =: r^{1,\infty} - J.$$

The first part of the rhs has already been estimated (prop. 2). It remains to estimate the last term $J$, that solves the following system :

$$\begin{cases} -\Delta J = 0, & \forall x \in \Omega^0, \\ J = 0, & \forall x \in \Gamma^0, \\ J = \epsilon \dfrac{\partial u^1}{\partial x_2}(x_1, 0)\left(\beta\left(\dfrac{x_1}{\epsilon}, \dfrac{1}{\epsilon}\right) - \overline{\beta}\right), & \forall x \in \Gamma^1, \\ J \text{ is } x_1 - \text{periodic on } \Gamma_{\text{in}} \cup \Gamma_{\text{out}}. \end{cases}$$

Using a $y_2$-linear lift $s$ that takes away the $\Gamma^1$ boundary term (which is exponentially small wrt $\epsilon$), and thanks to the Poincaré inequality, we obtain

$$\|J\|_{L^2(\Omega^0)} \leq c_{13}\|J\|_{H^1(\Omega^0)} \leq c_{14} e^{-\frac{1}{\epsilon}},$$

where $c_{13}$ and $c_{14}$ are constants independent on $\epsilon$. $\square$

REMARK 5.1. *The error in $O(\epsilon^{\frac{3}{2}})$ is only due to the first order boundary layer approximation. Indeed the extension of the Poiseuille flow is only linear inside $\Omega^\epsilon \setminus \Omega^0$. Nevertheless, we avoid errors when neglecting microscopic oscillations in our macroscopic problem as it was the case for $u^1$ and $u^2$.*

**5.2. A second order explicit wall-law.** Extending the same ideas as in the subsection above, one sets the following multi-scale problem: find $\mathcal{V}_\epsilon \in H^1(\Omega^0)$ such that

$$\begin{cases} -\Delta \mathcal{V}_\epsilon = C, & \forall x \in \Omega^0, \\ \mathcal{V}_\epsilon = \epsilon \dfrac{\partial u^2}{\partial x_2}\beta\left(\dfrac{x_1}{\epsilon}, 0\right) + \dfrac{\epsilon^2}{2}\dfrac{\partial^2 u^2}{\partial x_2^2}\gamma\left(\dfrac{x_1}{\epsilon}, 0\right), & \forall x \in \Gamma^0, \\ \mathcal{V}_\epsilon = 0, & \forall x \in \Gamma^1, \quad \mathcal{V}_\epsilon \text{ is } x_1 - \text{periodic on } \Gamma_{\text{in}} \cup \Gamma_{\text{out}}, \end{cases} \quad (5.2)$$

for which we can prove

PROPOSITION 7. *Under hypotheses 2.1, one gets*

$$\|u^\epsilon - \mathcal{V}_\epsilon\|_{L^2(\Omega^0)} \leq c_{15} e^{-\frac{1}{\epsilon}},$$

*where the constant $c_{15}$ is independent on $\epsilon$.*

**5.3. First order implicit wall-laws.** Note that the standard averaged wall-laws $u^1, u^2$ are building blocks of explicit multi-scale approximations $\mathcal{U}_\epsilon, \mathcal{V}_\epsilon$ solving problems (5.1,5.2). In this part we look for an implicit approximation that avoids the computation of these lower order approximations. Indeed, at first order we propose to solve :

$$\begin{cases} -\Delta \Upsilon_\epsilon = C, & \forall x \in \Omega^0, \\ \Upsilon_\epsilon = \epsilon\beta(\dfrac{x_1}{\epsilon}, 0)\dfrac{\partial \Upsilon_\epsilon}{\partial x_2}, & \forall x \in \Gamma^0, \\ \Upsilon_\epsilon = 0, & \forall x \in \Gamma^1, \quad \Upsilon_\epsilon \text{ is } x_1 - \text{periodic on } \Gamma_{\text{in}} \cup \Gamma_{\text{out}}. \end{cases} \quad (5.3)$$

We give here a first result of this kind :

THEOREM 5.1. *Under hypotheses 2.1, there exists a unique solution $\Upsilon_\epsilon \in H^1_{\Gamma^1}(\Omega^0)$ of problem (5.3). Moreover, one gets :*

$$\|u^\epsilon - \Upsilon_\epsilon\|_{L^2(\Omega^0)} \leq c_{16}\epsilon^{\frac{3}{2}}.$$

where $c_{16}$ is a constant independent of $\epsilon$.

*Proof.* There exists a unique solution $\Upsilon_\epsilon$ solving (5.3). Indeed, under hypotheses 2.1, the weak formulation of (5.3) reads :

$$a(u,v) := (\nabla u, \nabla v)_{\Omega^0} + \left(\frac{\partial u}{\partial x_2}, v\right)_{\Gamma^0}$$

$$= (\nabla u, \nabla v)_{\Omega^0} + \left(\frac{u}{\epsilon\beta}, v\right)_{\Gamma^0} = (C, v)_{\Omega^0} =: l(v), \quad \forall v \in H^1_{\Gamma^1}(\Omega^0),$$

At the microscopic level, we suppose that $P^0$ does not cross $\Gamma$, thus there exists a minimal distance $\delta > 0$ separating them. By the maximum principle, $\beta$ is bounded: $\beta \in [\delta; 1]$. Thus $1/\beta$ is bounded a.e. The bilinear form $a$ is continuous coercive in $H^1_{\Gamma^1}(\Omega^0)$, the linear form $l$ is continuous as well, thus existence and uniqueness follow by the Lax-Milgram theorem. To estimate this new approximation's convergence rate we add and substract $\mathcal{U}_\epsilon$, the explicit wall-law between $u^\epsilon$ and $\Upsilon_\epsilon$.

$$r^1_{\text{bl},i} := u^\epsilon - \Upsilon_\epsilon = u^\epsilon - \mathcal{U}_\epsilon + \mathcal{U}_\epsilon - \Upsilon_\epsilon = r^1_{\text{bl}} + \mathcal{U}_\epsilon - \Upsilon_\epsilon =: r^1_{\text{bl}} + \Theta. \quad (5.4)$$

$\Theta$ is the solution of the boundary value problem reading :

$$\begin{cases} -\Delta\Theta = 0, & \forall x \in \Omega, \\ \Theta = \epsilon\beta\left[\dfrac{\partial u^1}{\partial x_2} - \dfrac{\partial \Upsilon_\epsilon}{\partial x_2}\right], & \forall x \in \Gamma^0, \\ \Theta = 0, & \forall x \in \Gamma^1, \quad \Theta \text{ is } x_1 - \text{periodic on } \Gamma_{\text{in}} \cup \Gamma_{\text{out}}. \end{cases}$$

We reexpress the boundary condition on $\Gamma^0$ introducing a Robin like condition, namely :

$$\Theta - \epsilon\beta\frac{\partial\Theta}{\partial x_2} = \epsilon\beta\left[\frac{\partial u^1}{\partial x_2} - \frac{\partial \mathcal{U}_\epsilon}{\partial x_2}\right], \quad \forall x \in \Gamma^0, \quad (5.5)$$

where the rhs is explicitly known. We have the following weak formulation :

$$-(\Delta\Theta, v)_{\Omega^0} = -\left(\frac{\partial\Theta}{\partial \mathbf{n}}, v\right)_{\Gamma^0} + (\nabla\Theta, \nabla v)_{\Omega^0} = 0, \quad \forall v \in H^1_{\Gamma^1}(\Omega^0),$$

where the space $H^1_{\Gamma^1}(\Omega^0)$ contains $H^1(\Omega^0)$ functions vanishing on $\Gamma^1$. Then using (5.5) one writes :

$$a(\Theta, v) = (\nabla\Theta, \nabla v)_{\Omega^0} + \left(\frac{\Theta}{\epsilon\beta}, v\right)_{\Gamma^0} = \left(\frac{\partial u^1}{\partial x_2} - \frac{\partial \mathcal{U}_\epsilon}{\partial x_2}, v\right)_{\Gamma^0}.$$

We remark that the rhs is in fact a boundary term of another comparison problem and we set $z = u^1 - \mathcal{U}_\epsilon$ where $z$ is harmonic and solves :

$$\left(\frac{\partial z}{\partial x_2}, v\right)_{\Gamma^0} = -(\Delta z, v)_{\Omega^0} - (\nabla z, \nabla v)_{\Omega^0}, \quad \forall v \in H^1_{\Gamma^1}(\Omega^0).$$

*Estimates of the gradient.* We have recovered a simpler problem that reads

$$a(\Theta, v) = -(\nabla z, \nabla v)_{\Omega^0}, \quad \forall v \in H^1_{\Gamma^1}(\Omega^0).$$

Thanks to proposition 2 and proposition 6, one gets

$$\|\nabla\Theta\|_{L^2(\Omega^0)} \leq \|\nabla z\|_{L^2(\Omega^0)} \leq \|\nabla(u^\epsilon - u^1)\|_{L^2(\Omega^\epsilon)} + \|\nabla(u^\epsilon - \mathcal{U}_\epsilon)\|_{L^2(\Omega^0)} \leq 2c_{17}\epsilon,$$

where $K$ is a constant independent of $\epsilon$.

*Estimate of the trace.* The control on the interior term enables to recover trace estimates

$$\|\Theta\|_{L^2(\Gamma^0)}^2 \leq \|\beta\|_{L^\infty(\Gamma)} \int_0^L \frac{\Theta^2(x_1,0)}{\beta\left(\frac{x_1}{\epsilon},0\right)} dx_1 \leq \epsilon \|\nabla\Theta\|_{L^2(\Omega^0)} \|\nabla z\|_{L^2(\Omega^0)} \leq c_{17}^2 \epsilon^3.$$

*Final estimate.* By the dual problem, and trace estimates above, we finaly obtain

$$\|\Theta\|_{L^2(\Omega^0)} \leq c_{18} \|\Theta\|_{L^2(\Gamma^0)} \leq c_{19} \epsilon^{\frac{3}{2}},$$

Recalling relation (5.4), one gets :

$$\|r^1_{\text{bl},i}\|_{L^2(\Omega^0)} \leq \|r^1_{\text{bl}}\|_{L^2(\Omega^0)} + \|\Theta\|_{L^2(\Omega^0)},$$

which ends the proof. □

REMARK 5.2. *A similar implicit approach could be considered at second order. This should lead to consider a multi-scale Wentzel condition. It is an open problem to show existence, uniqueness and error estimates as in theorem 5.1 in this case.*

**6. Numerical evidence.** We compute $u_\Delta^\epsilon$, a numerical approximation of the rough problem (2.2) on the whole domain $\Omega^\epsilon$, $\epsilon$ taking a given range of values in $[0.1, 1]$. Then, we restrict the computational domain to $\Omega^0$, and compute macroscopic approximations $u_\Delta^1, u_\Delta^2, \mathcal{U}_{\epsilon,\Delta}, \mathcal{V}_{\epsilon,\Delta}, \Upsilon_{\epsilon,\Delta}$, again for each value of $\epsilon$. We evaluate the errors w.r.t. $u_\Delta^\epsilon$ interpolating the latter exact solution over the meshes of the former ones.

*Computational setting.* For every simulation, we use a $\mathbb{P}_2$ Lagrange finite element method implemented in the C++ code rheolef[1] [21]. Our computational domain is a channel of length $L = 10$ and of height $h = 1$. We assume a rough periodic bottom boundary $\Gamma^\epsilon$ defined by formula (2.1) with

$$f(y_1) := -\frac{(1 + \cos(y_1))}{2} - \delta,$$

where $\delta$ is a positive constant set to $5e - 2$.

*The rough solution $u_\Delta^\epsilon$.* We compute $u_\Delta^\epsilon$ over a single macroscopic cell $x \in \omega^\epsilon := \{x_1 \in [0, 2\pi\epsilon] \text{ and } x_2 \in [f(x_1/\epsilon), 1]\}$ and we assume periodic boundary conditions at $\{x_2 = 0\} \cup \{x_1 = 2\pi\epsilon\}$. For each fixed $\epsilon$, we mesh the domain $\omega^\epsilon$ while keeping approximately the same number of vertices in the $x_1$ direction. This forces the mesh to get finer in the $x_2$ direction in order to preserve the ratio between the inner and outer radius of each triangular element. With such a technique we avoid discretizations that could be of the same order as $\epsilon$.

*Cell problems.* In order to extract fruitful information for macroscopic wall-laws, we compute first and second order cell problems. Again we impose $y_1$-periodic boundary conditions. We truncate the upper infinite part of the domain by imposing a homogeneous Neumann boundary condition at $y_2 = 10$ after verifying that a variation of the domains height no more affects the results. In [14], the authors show an exponential convergence w.r.t. to the height of the truncated upper domain towards the $y_2$-infinite $y_1$-periodic cell problems (3.2) and (3.6), this validates our approach. Cell problems are computed over a mesh containing (9211 elements and 4738 vertices).

---
[1] http://ljk.imag.fr/membres/Pierre.Saramito/rheolef/

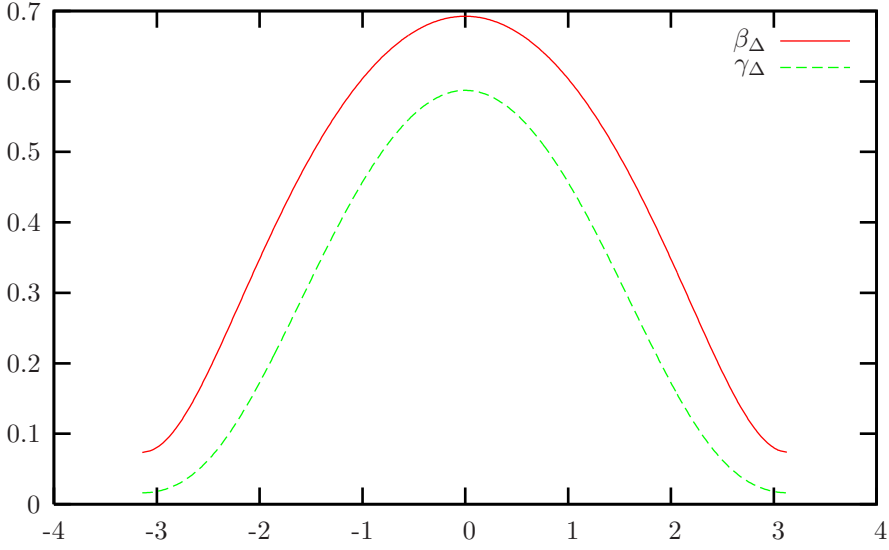

figure 6.1: The traces $\beta_\Delta(y_1,0)$ and $\gamma_\Delta(y_1,0)$

| $e_\Delta$ | $u_\Delta^\epsilon - u_\Delta^0$ | $u_\Delta^\epsilon - u_\Delta^1$ | $u_\Delta^\epsilon - u_\Delta^2$ | $u_\Delta^\epsilon - \mathcal{U}_{\epsilon,\Delta}$ | $u_\Delta^\epsilon - \mathcal{V}_{\epsilon,\Delta}$ | $u_\Delta^\epsilon - \Upsilon_{\epsilon,\Delta}$ |
|---|---|---|---|---|---|---|
| $\alpha$ | 1.11 | 1.4786 | 1.3931 | 1.768 | 2-3.6 | 1.6227 |

TABLE 6.1
*Numerical orders of convergence for various approximations*

We extract solutions' trace on the fictitious interface $\Gamma$ for both first and second order cell problems (cf. fig. 6.1), and compute the averages $\overline{\beta} = 0.43215$ and $\overline{\gamma} = 0.29795$.
*Macroscopic approximations: Classical & new wall-laws* We compute the classical macroscopic wall-laws over $\omega_+^\epsilon = \{x \in \omega^\epsilon \, / \, x_2 \geq 0\}$, a single periodicity cell of $\Omega^0$. We follow the same rate of refinement as described above. Then, we solve problems (4.2,4.4).
In the same spirit, we use both averages $(\overline{\beta}, \overline{\gamma})$ and oscillating functions $\beta(\frac{x_1}{\epsilon},0), \gamma(\frac{x_1}{\epsilon},0)$ as a non-homogenous Dirichlet boundary condition over the macroscopic domain when solving (5.1) and (5.2). To provide values at the boundary we use a $\mathbb{P}_1$ interpolation of the data extracted from the cell problems.
For the implicit multi-scale wall-law, we solve system (5.3) using the inverse of $\beta_\Delta(x_1/\epsilon,0)$ as a weight in the boundary integrals of the discrete variational formulation.

*Results.* We plot fig. 6.2, the $L^2(\Omega^0)$ error computed respectively for approximations presented above: $u_\Delta^\epsilon - u_\Delta^0, u_\Delta^\epsilon - u_\Delta^1, u_\Delta^\epsilon - u_\Delta^2, u_\Delta^\epsilon - \mathcal{U}_{\epsilon,\Delta}, u_\Delta^\epsilon - \mathcal{V}_{\epsilon,\Delta}, u_\Delta^\epsilon - \Upsilon_{\epsilon,\Delta}$. If we set $e_\Delta = C\epsilon^\alpha$, table 6.1 gives approximate numeric values of convergence rates.
*Interpretation.* A first important result, visible fig. 6.2, is that there is no difference between first and second order macroscopic wall-laws $u^1$ and $u^2$. This proves that our estimates are actually optimal. It explains also why one could never distinguish first from second order approximations in [2, 1].
Next, we remark that convergence orders are not better than those predicted by the estimates for $u_\Delta^1, u_\Delta^2, \mathcal{U}_{\epsilon,\Delta}$, while the error displayed for $\mathcal{V}_{\epsilon,\Delta}$ is limited by the $\mathbb{P}_2$ interpolation. Indeed, the $H^1(\Omega^0)$ error is of order 3 on the vertices but is worse elsewhere inside the elements. Nevertheless, the error $u_\Delta^\epsilon - \mathcal{V}_{\epsilon,\Delta}$ is more than one order smaller than for $u_\Delta^1, u_\Delta^2, \mathcal{U}_{\epsilon,\Delta}$ for every fixed $\epsilon$.
The fully explicit oscillating wall-laws $\mathcal{U}_{\epsilon,\Delta}, \mathcal{V}_{\epsilon,\Delta}$ provide better results than the implicit ones, $u^1, u^2$ and $\Upsilon_\epsilon$. Indeed, in the former the shear rate $\partial u^0/\partial x_2(x_1,0)$ and the second order derivative $\partial^2 u^0/\partial x_2^2(x_1,0)$ of the limit Poiseuille profile are explicit and included in the boundary condition, whereas the latter approximate this information as well. This leads to supplementary errors on the macroscopic scale for implicit wall-laws.

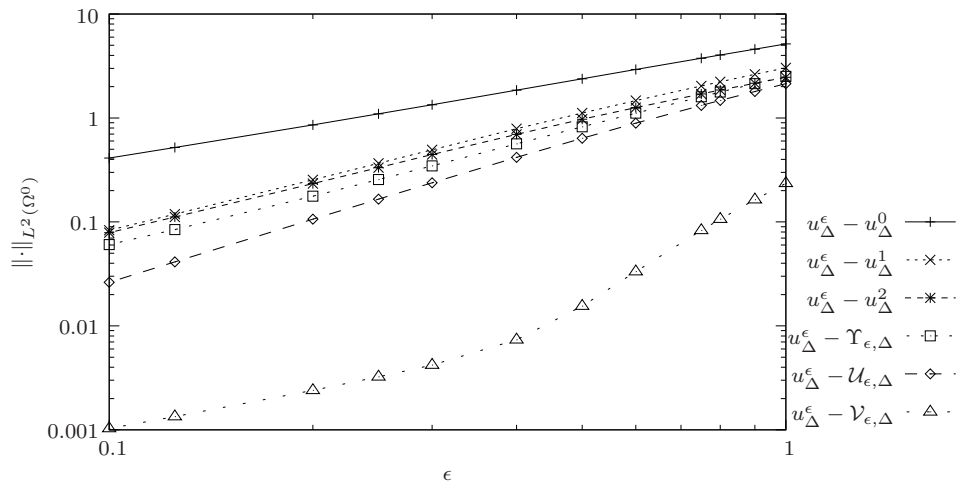

figure 6.2: $L^2(\Omega^0)$ error computed versus $\epsilon$

**Acknowledgments.** The first author is partially supported by the project "Études mathématiques de paramétrisations en océanographie" that is part of the "ACI jeunes chercheurs 2004" framework of the French Research Ministry and by a Rhône Alpes project "Equations de type Saint-Venant avec viscosité pour des problèmes environnementaux". The second author was partially supported by a contract with Cardiatis® a company providing metallic multi-layer stents for cerebral and aortic aneurysms. This research has been partly funded by the Rhône-Alpes Institute of Complex Systems IXXI[2]. The authors would like to thank E. Bonnetier for fruitful discussions and helpful proofreading.

---

[2]http://www.ixxi.fr/

## Appendix A. The cell problems.

**A.1. Various properties of the first order cell problem's solution.** Existence and uniqueness of solutions of system (3.2), have been partially proven in [2]. The authors consider a truncated domain supplied with a non-local "transparency" condition, the latter is obtained via the fourier transform. We give here a rigorous proof in the unbounded domain framework.

*Proof.* [of theorem 3.1] In what follows we express the cell problem as an inverse Steklov-Poincaré problem solved on the fictitious interface $\Gamma$. This allows us to characterize $\beta$ the solution of (3.2) on domains $Z^+$ and $P$ separately, as depending only on $\eta$, the trace on $\Gamma$. We apply domain decomposition techniques [20]. In a first step we give a simple proof of existence that guarantees the existence of the gradient in $L^2(Z^+ \cup P)$. The solutions of the cell problems are not in the classical Sobolev spaces because the domain $Z^+$ is unbounded in the $y_2$ direction: the solutions are only locally integrable. For this purpose, we define, for an arbitrary open set $\omega$,

$$D^{n,p}(\omega) = \{v \in L^1_{\text{loc}}(\omega) / \, D^\alpha v \in \mathrm{L}^p(\omega), \, \forall \alpha \in \mathbb{Z}^{d^n} |\alpha| = n, \, v \text{ is } y_1 - \text{periodic }\}.$$

In the particular case when $n = 1$ and $p = 2$, we define $D^{1,2}_0(\omega) := \{v \in D^{1,2}(\omega) / \, v|_{\partial\omega} = 0\} := V_0(\omega)$, which is a Hilbert space for the norm of the gradient.

LEMMA A.1. *Problem (3.2) admits a unique solution $\beta$ belonging to $D^{1,2}(Z^+ \cup P)$.*

*Proof.* [of lemma A.1] We define the lift $s = y_2 \chi_{[P]}$ that belongs to $D^{1,2}(Z^+ \cup P)$. Setting $\tilde{\beta} = \beta - s$, the lifted problem becomes

$$\begin{cases} -\Delta\tilde{\beta} = \delta_\Gamma, & \text{in } Z^+ \cup P, \\ \tilde{\beta} = 0, & \text{on } P^0, \end{cases}$$

where $\delta_\Gamma$ is the dirac mesure that concentrates on the fictitious interface $\Gamma$. The equivalent variational form of this problem reads

$$a(\tilde{\beta}, v) = l(v), \quad v \in D_0^{1,2}(Z^+ \cup P), \tag{A.1}$$

where $a(u,v) = (\nabla u, \nabla v)_{Z^+ \cup P}$ and $l(v) = -(\nabla s, \nabla v)_P$. These forms are obviously continuous bilinear (resp. linear) on $D_0^{1,2}(Z^+ \cup P) \times D_0^{1,2}(Z^+ \cup P)$ (resp. $D_0^{1,2}(Z^+ \cup P)$). Because of the homogeneous boundary condition the semi-norm of the gradient is a norm. By Lax-Milgram theorem, the desired result follows. □

We define the following spaces

$$V_1 = D^{1,2}(Z^+), \quad V_2 = \{v \in H^1(P) \text{ s.t. } v|_{P^0} = 0, \quad v \text{ is } y_1-\text{periodic }\}$$
$$V_{1,0} = \{v \in V_1, \quad v|_\Gamma = 0\}, \quad V_{2,0} = \{v \in V_2, \quad v|_\Gamma = 0\}$$
$$\Lambda = \{\eta \in H^{\frac{1}{2}} \text{ s.t. } \eta = v|_\Gamma \text{ for a suitable } v \in D_0^{1,2}(Z^+ \cup P)\}.$$

LEMMA A.2. *The following domain decomposition problem is equivalent to* (A.1) : *we look for* $(\beta_1, \beta_2) \in V_1 \times V_2$ *such that*

$$\begin{cases} a_1(\beta_1, v) := (\nabla\beta_1, \nabla v)_{Z^+} = 0, & \forall v \in V_{1,0}, \\ \beta_1 = \beta_2, & \text{on } \Gamma, \\ a_2(\beta_2, v) := (\nabla\beta_1, \nabla v)_P = -(\nabla s, \nabla v)_P \equiv 0, & \forall v \in V_{2,0}, \\ a_2(\beta_2, \mathcal{R}_2\mu) = -(\nabla s, \nabla\mathcal{R}_2\mu) - a_1(\beta_1, \mathcal{R}_1\mu), & \forall \mu \in \Lambda, \end{cases} \tag{A.2}$$

*where $\mathcal{R}_i$ denotes any possible extension operator from $\Gamma$ to $V_i$.*

*Proof.* [of lemma A.2] Let us start by considering the solution $\beta$ of (A.1). Setting $\beta_1 = \beta|_{Z^+}, \beta_2 = \beta|_P$, we have that $\beta_i \in V_i$ and that (A.2).1,(A.2).2 and (A.2).3 are trivially satisfied. Moreover, for each $\mu \in \Lambda$, the function $\mathcal{R}\mu$ defined as $\mathcal{R}\mu = \mathcal{R}_1\mu\chi_{Z^+} + \mathcal{R}_2\mu\chi_P$ belongs to $V_0$. Therefore we have $a(\beta, \mathcal{R}\mu) = (f, \mathcal{R}\mu), \forall \mu \in \Lambda$ which is equivalent to (A.2).4.

On the other hand, let $\beta_i$ be the solution of (A.2). Setting $\beta = \beta_1\chi_{[Z^+]} + \beta_2\chi_{[P]}$ from (A.2).2, it follows that $\nabla\beta \in L^2(Z^+ \cup P)$, and $\beta|_{P^0} = 0$. Then taking $v \in V_0$ we set $\mu = v|_\Gamma \in \Lambda$. Define $\mathcal{R}\mu$ as before; clearly $(v_i - \mathcal{R}_i\mu) \in V_{i,0}$ and from (A.2).1, (A.2).3, (A.2).4 it follows that

$$a(\beta, v) = \sum_i [a_i(\beta_i, v_i - \mathcal{R}_i\mu) + a_i(\beta_i, \mathcal{R}_i\mu)] = -(\nabla s, \nabla\mathcal{R}_2\mu)_P$$
$$= -(1, \mu)_\Gamma = -(1, v)_\Gamma = -(\nabla s, \nabla v)_P.$$

□

*The Steklov-Poincaré operator.* The Steklov-Poincaré operator $S$ acts between the space of trace functions $\Lambda$ and its dual. More precisely, applying Green's formula and setting $H_i\eta$ to be the harmonic lift in $Z^+$ (resp. $P$) for all $\eta \in \Lambda$, we have

$$< S\eta, \mu > = \sum_i \left\langle \frac{\partial}{\partial \nu_i} H_i\eta, \mu \right\rangle = \int_{Z^+} \nabla H_1\eta \cdot \nabla \mathcal{R}_1\mu + \int_P \nabla H_2\eta \cdot \nabla \mathcal{R}_2\mu$$
$$= \sum_i a_i(H_i\eta, \mathcal{R}_i\mu), \quad \forall \eta, \mu \in \Lambda,$$

where $< \cdot, \cdot >$ denotes the duality pairing between $\Lambda'$ and $\Lambda$. In particular, taking $\mathcal{R}_i\mu = H_i\mu$, we obtain the following variational representation:

$$< S\eta, \mu > = \sum_i a_i(H_i\eta, H_i\mu), \quad \forall \eta, \mu \in \Lambda.$$

*The linear form on $\Lambda$.* We set $l(\mu)$ as follows:

$$l(\mu) = -(\nabla s, \nabla H_2\mu)_P = \left(1, \frac{\partial}{\partial x_2} H_2\mu\right)_P = (1, \mu)_\Gamma.$$

LEMMA A.3. *The problem:*

$$\text{find } \eta \in \Lambda \text{ such that } < S\eta, \mu > = l(\mu), \quad \forall \mu \in \Lambda, \tag{A.3}$$

*admits a unique solution. Moreover this is equivalent to solve* (A.2).

*Proof.* [of lemma A.3] We use the Lax-Milgram framework:
- Continuity:

$$< S\eta, \mu > \leq \|\nabla H\eta\|_{L^2(Z^+\cup P)} \|\nabla H\mu\|_{L^2(Z^+\cup P)} \leq c_{20} \|\eta\|_\Lambda \|\mu\|_\Lambda,$$

by well know estimates for solutions of elliptic boundary value problems [15]. For $H_2$ this can be computed explicitly (see below). The continuity of $l$ is obvious.
- Coercivity

$$< S\eta, \eta > = \|\nabla H\eta\|^2_{L^2(Z^+\cup P)} \geq c_{21} \|H_2\eta\|^2_{H^1(P)} \geq c_{22} \|\eta\|^2_\Lambda.$$

Then applying Lax-Milgram theorem one gets the desired result.
To prove the equivalence between (A.3) and (A.2), it suffices to separate the harmonic lift $H_i$ and the solutions of the Poisson problem with homogeneous boundary conditions and the result follows as in [20] p.10. ∎

The harmonic extension in $Z^+$ named $H_1$ We set $\eta \in \Lambda$. By decomposing in $y_1$-fourier modes, one gets that the solution of:

$$\begin{cases} \Delta \beta = 0, & \forall y \in Z^+, \\ \beta = \eta, & \forall y \in \Gamma, \end{cases} \tag{A.4}$$

rewritten as $\beta = \sum_k \beta_k(y_2) e^{iky_1}$, $\forall y \in Z^+$ should satisfy the following system of ODE's:

$$\begin{cases} \beta_k'' - k^2 \beta_k = 0, & y_2 \in \mathbb{R}^+ \\ \beta_k(0) = \eta_k, & y_2 = 0 \\ \beta_k(y_2) \in L^\infty(\mathbb{R}^+; \mathbb{C}), \end{cases}$$

where $\eta_k = \int_0^{2\pi} e^{-iky_1}\eta(y_1)dy_1$ are $\eta$'s fourier coefficients on $\Gamma$. The solution $\beta_{Z^+}$ is explicit and reads

$$H_1\eta = \beta_{|Z^+} = \sum_{k=-\infty}^{\infty} \eta_k e^{-|k|y_2+iky_1}, \quad \forall y \in Z^+. \tag{A.5}$$

To show exponential convergence towards zero of $\beta - \overline{\beta}$ and $\nabla\beta$ when $y_2 \to 0$, we use the same arguments as in the second part of [3], theorem 2.2.1 p. 637, whose proof is omitted.

PROPOSITION 8. *There exists $\alpha_1 \geq (4\pi)^2/9$ such that the solution of problem (3.2) satisfies*

$$\left\|\beta - \overline{\beta}\right\|_{L^2(Z^+\cup P, e^{\alpha_1 y_2})} \leq c_{23}\|\nabla\beta\|_{L^2(Z^+\cup P, e^{\alpha_1 y_2})} \leq c_{24},$$

*which implies also $\beta$'s and $\nabla\beta$'s exponential decay in the $y_2$ direction.*
☐

**A.2. The second order boundary layer.** *Proof.* [of proposition 3] Problem (3.6) is equivalent to solve :

$$\begin{cases} \Delta\tilde{\gamma} = 2\chi_{[P]}, & \forall y \in Z^+ \cup P, \\ \tilde{\gamma} = 0, & \forall y \in P^0. \end{cases}$$

This, under the previous domain decomposition form, reads: find $(\tilde{\gamma}_{Z^+}, \tilde{\gamma}_P)$ such that

$$\begin{cases} (\nabla\tilde{\gamma}_{Z^+}, \nabla v)_{Z^+} = 0, & \forall v \in H^1_\Gamma(Z^+), \\ \tilde{\gamma}_{Z^+} = \tilde{\gamma}_P, & \text{on } \Gamma, \\ (\nabla\tilde{\gamma}_P, \nabla v)_P = -(2, v)_P, & \forall v \in H^1_{\Gamma\cup P^0}(P), \\ (\nabla\tilde{\gamma}_P, \nabla\mathcal{R}_P\mu)_P = -(2, \mathcal{R}_P\mu)_P - (\nabla\tilde{\gamma}_{Z^+}, \mathcal{R}_{Z^+}\mu)_{Z^+}, & \forall \mu \in H^{\frac{1}{2}}(\Gamma). \end{cases} \tag{A.6}$$

Following the same lines as the proof above, we write the interface problem :

$$< S\lambda, \mu > = (\nabla H_P\lambda, \nabla H_P\mu) + (\nabla H_{Z^+}\lambda, \nabla H_{Z^+}\mu), \quad \forall \mu \in H^{\frac{1}{2}}(\Gamma),$$
$$= -(2, H_P\mu) - (\nabla\mathcal{G}_2, \nabla H_P\mu) =: l(\mu), \quad \forall \mu \in H^{\frac{1}{2}}(\Gamma),$$

where $\mathcal{G}_2$ is the solution of the homogeneous Poisson problem :

$$\begin{cases} \Delta\mathcal{G}_2 = 2, & \forall y \in P, \\ \mathcal{G}_2 = 0, & \forall y \in P^0 \cup \Gamma, \\ \mathcal{G}_2 & \text{is } y_1 - \text{periodic}. \end{cases}$$

One gets the continuity of the linear form again, thanks to the properties of the harmonic lifts [15, 8] :

$$|l(\mu)| = |-(2, H_P\mu) - (\nabla\mathcal{G}_2, \nabla H_P\mu)| \leq c_{25}\|H_P\mu\|_{H^1(P)} \leq c_{26}\|\mu\|_{H^{\frac{1}{2}}(\Gamma)}.$$

And again, by the Lax-Milgram theorem, one gets the desired result. ☐